\documentclass{amsart}
\usepackage{graphicx}
\newtheorem{theorem}{Theorem}

\newtheorem{corollary}[theorem]{Corollary}

\newtheorem{lemma}[theorem]{Lemma}

\newcommand{\abs}[1]{\left\vert#1\right\vert}

\usepackage{amsmath,amscd}

\begin{document}
\title[]{Complex Tangencies to Embeddings of Heisenberg Groups and Odd-Dimensional Spheres}%
\author{Ali M. Elgindi}%
%\thanks{}%
%\subjclass{}%
%\keywords{}%

%\date{}%
%\dedicatory{}%
%\commby{}%
% ----------------------------------------------------------------
\begin{abstract}
 The notion of a complex tangent arises for embeddings of real manifolds into complex spaces. It is of particular interest when studying embeddings of real $n$-dimensional manifolds into $\mathbb{C}^n$. The generic topological structure of the set complex tangents to such embeddings $M^n \hookrightarrow \mathbb{C}^n$ takes the form of a (stratified) $(n-2)$-dimensional submanifiold of $M^n$. In this paper, we generalize our results from our previous work for the 3-dimensional sphere and the Heisenberg group to obtain results regarding the possible topological configurations of the sets of complex tangents to embeddings of odd-dimensional spheres $S^{2n-1} \hookrightarrow \mathbb{C}^{2n-1}$ by first considering the situation for the higher dimensional analogues of the Heisenberg group.
\end{abstract}
\maketitle

\par\ \par\

% --------------------------------------------------------------

\section*{0. Introduction}
The study of CR-manifolds and the local analytic structure of real submanifolds of complex space has been a subject of extensive study over last half-century, originated and inspired by E. Bishop's pioneering paper in the 1960's (see [1]). Of particular interest for many mathematicians was the situation in dimension 2, which is now well-understood via results such as those derived by F. Forstneric in his paper [3], and other works. The structure of complex tangents to embeddings of higher dimensional manifolds is not as well studied, although there have been some interesting results in this field as well. In particular one could see the work of H.F. Lai in his 1972 paper (cited [5]), where he derived formulae relating the characteristic classes of a manifold immersed into a complex manifold with the fundamental class of the set of complex tangents to the immersion. X. Gong also investigated the existence of totally real spheres of real dimension less than $n$ into the complex Euclidean space $\mathbb{C}^n$ in his paper ([4]). Our investigations are along a different direction than those of Lai and Gong, in particular we are interested in the "critical" case where the dimension of the real manifold being embedded is the same as the complex dimension of the space into which it is being embedded.
\par\ \par\
In our previous work (cited [2]), we proved that every topological type of knot in $S^3$ may arise as the set of complex tangents to an embedding $S^3 \hookrightarrow \mathbb{C}^3$. We accomplished this by first proving the result in the Heisenberg group $\mathbb{H}$ using its tangential Cauchy-Riemann operator. We then used the natural sterographic projection of $S^3$ and a certain biholomorphic map to obtain the result for $S^3$.
\par\ \par\
The situation in higher dimensions is more complicated, however. There are several Cauchy-Riemann operators in higher dimensions, and to analyze the set of complex tangents we will need to consider the zero sets of each one.
\par\ \par\
For the higher dimensional analogues of the Heisenberg group, which we denote by $\mathbb{H}^{2n-1}$ (for $n \geq 3$), the situation may be simplified to consider only one of the Cauchy-Riemann operators. By using our results for $\mathbb{H} = \mathbb{H}^3$ (from [2]), we are able to prove that every real algebraic subset given by one or two polynomial equations in $\mathbb{H}^{2n-1}$ may arise as the set of complex tangents to some embedding of $\mathbb{H}^{2n-1}$ into $\mathbb{C}^{2n-1}$.
\par\ \par\
We note that the Heisenberg groups are naturally equivalent to real Euclidean spaces, so our above result merely asserts that we can construct embeddings of real Euclidean space which are complex tangent along a given subset. However, we will make use of the natural stereographic projection of $S^{2n-1}$ onto $\mathbb{H}^{2n-1}$ to prove the analogous result for all odd-dimensional spheres, in particular that every real algebraic subset given by one or two polynomial equations in $S^{2n-1}$ may arise as the set of complex tangents to some embedding of $S^{2n-1}$ into $\mathbb{C}^{2n-1}$. We also obtain a result regarding the "minimal" configuration of complex tangents to embeddings of spheres.
\par\ \par\
We also note that we can generalize our results for the Heisenberg groups to "trough-like" hypersurfaces of $\mathbb{R}^{2n}$, and in particular we obtain totally real embeddings into $\mathbb{C}^{2n-1}$ for such hypersurfaces.
\par\ \par\
\section*{1. Higher Dimensional Heisenberg Groups}

Let $M^{2n-1} \hookrightarrow \mathbb{C}^n$ be an embedding of a real $(2n-1)$-dimensional manifold into complex $n$-space. The notion of tangential Cauchy-Riemann operators as we defined for hypersurfaces $M^3 \subset \mathbb{C}^2$ in [2] will admit a generalization to higher dimensions. In fact for $n > 2$ we expect "several" tangential CR-operators at any given point of $M$.
\par\ \par\
Let $M =\{\rho = 0\} \subset \mathbb{C}^n$, $\overline{\rho} = \rho$, $\nabla \rho \neq 0$ on $M$. A Cauchy-Riemann operator will be an operator of the form: $L_\textbf{a} = \Sigma_{j=1}^n \overline{a_j} (z) \frac{\partial}{\partial \overline{z_j}}$, where $\textbf{a} = (a_j): M \rightarrow \mathbb{C}^n$ is continuous and $z =(z_1,..., z_n)$ are the standard  (holomorphic) coordinates of $\mathbb{C}^n$. We say such a CR-operator $L$ is tangential to $M$ if $L(\rho) (\zeta) = 0$  $\forall \zeta \in M$, i.e.:
\par\
$\Sigma \overline{a_j} (\zeta) \frac{\partial \rho}{\partial \overline{z_j}} (\zeta) = 0$
\par\ \par\
We say a function $u: M \hookrightarrow \mathbb{C}$ satisfies the Cauchy-Riemann equations on $M$ if $Lu = 0$ for all tangential CR-operators to $M$. This condition is equivalent to the vanishing of the 2-form: $\overline{\partial} u \wedge \overline{\partial} \rho = 0$ over $M$; see Rudin's book [5] for reference.
\par\ \par\
Let us now turn our focus to the Heisenberg groups; to avoid confusion, we henceforth denote the "ordinary" Heisenberg group $\mathbb{H} \subset \mathbb{C}^2$ as $\mathbb{H}^3$. Let $M =\mathbb{H}^{2n-1}= \{Im(z_n) = \Sigma_{j=1}^{n-1} \abs{z_j}^2 \} \subset \mathbb{C}^n$.  We may then write $\rho = i(\overline{z_n} - z_n) - 2 (\Sigma_{j=1}^{n-1} z_j \overline{z_j})$.
\par\ \par\
Define the operators: $L_{jn} = \frac{\partial \rho}{\partial \overline{z_n}} \frac{\partial}{\partial \overline{z_j}} - \frac{\partial \rho}{\partial \overline{z_j}} \frac{\partial}{\partial \overline{z_n}}$. Note $L_{jn}$ corresponds to the vector function:
\par\
$\textbf{a}_{jn} = \frac{\partial \rho}{\partial z_n} \textbf{e}_j - \frac{\partial \rho}{\partial z_j} \textbf{e}_n$ by the above correspondence.
In fact, these operators $\{L_{1n},...,L_{n-1,n}\}$ form a basis for the CR-operators. Hence, a function u satisfies $\overline{\partial} u \wedge \overline{\partial} \rho = 0$  if and only if  $L_{jn} (u) = 0$ for every $j=1, ..., n-1$. Such functions $u$ are called CR functions.
\par\ \par\
Hence, a function $u: \mathbb{H}^{2n-1} \hookrightarrow \mathbb{C}$ satisfies the Cauchy-Riemann equations on $\mathbb{H}^{2n-1}$ if and only if  $L_{jn} (u) = 0$ for every $j=1, ..., n-1$. These $(n-1)$ differential operators then serve as the "relevant" tangential CR-operators to $\mathbb{H}^{2n-1}$.
\par\ \par\
We compute these operators explicitly as: $L_{jn} = 2z_j \frac{\partial}{\partial \overline{z_n}} + i \frac{\partial}{\partial \overline{z_j}}$. Using our proof for Lemma 4 in [2], we assert analogously:
\par\
\begin{lemma}: Acting on the space of polynomials $\mathcal{P}$ in the variables $\{z_k, \overline{z_k}\}_{k=1}^n$, the linear operators $L_{jn}:\mathcal{P} \rightarrow \mathcal{P}$ are onto, for every $j=1,...,n-1$.  \end{lemma}
\par\
$\textbf{\emph{\underline{Proof:}}}$
\par\
Fix $1 \leq j \leq n-1$. Since $L_{jn}$ is linear, it suffices to show that every monomial ${z_1}^{k_1} \overline{z_1}^{l_1} ... {z_n}^{k_n} \overline{z_n}^{l_n}$ is in $Im(L_{jn})$.
\par\ \par\
Furthermore, as $L_{jn}$ only acts on the variables $z_j, \overline{z_j}, z_n, \overline{z_n}$, it suffices to show that every term: ${z_j}^{k_1} \overline{z_j}^{l_1} {z_n}^{k_2} \overline{z_n}^{l_2}$ is in the range of $L_{jn}$. But the proof of this fact will follow in a similar fashion as the proof of the analogous result in two complex variables, which we gave in Lemma 4 in our paper [2].
\par\ \par\
For the proof, let us first note that as $L_{jn} ({z_n}^r g(z_j,z_n))= {z_n}^r L_{jn}(g(z_j,z_n))$, we need only verify the identity for all monomials: $f(z_j,z_n)={z_j}^m \overline{z_j}^k \overline{z_n}^l$.
\par\
We proceed by induction on $ l \in \mathbb{N}$. Note if $l=0$, then $f ={z_j}^m \overline{z_j}^k = L_{jn}
(\frac{1}{i(k+1)} {z_j}^m \overline{z_j}^{k+1})$, for any $m, k \in \mathbb{N}$.
Now, let $l \in \mathbb{N}$ be fixed and suppose $z_j^m \overline{z_j}^k \overline{z_n}^r \in L_{jn}(\mathcal{P})$, for all $m, k \in \mathbb{N}$, and $r\leq l$.
We wish to show $f(z_j,z_n)=z_j^m \overline{z_j}^k \overline{z_n}^{l+1} \in L_{jn}(\mathcal{P})$. We have that
\par\
$L_{jn}(\frac{1}{i(k+1)} z_j^m \overline{z_j}^{k+1} \overline{z_n}^{l+1}) = z_j^m \overline{z_j}^k \overline{z_n}^{l+1} + \frac{2(l+1)}{i(k+1)} z_j^{m+1} \overline{z_j}^{k+1} \overline{z_n}^l$.
\par\
Note that the second term on the right-hand side is in the range of $L_{jn}$ by the induction hypothesis and linearity. Let $g: \mathbb{C}^n \rightarrow \mathbb{C}$ be a polynomial such that $L_{jn} (g)= \frac{2(l+1)}{i(k+1)} z_j^{m+1} \overline{z_j}^{k+1} \overline{z_n}^l$.
The term on the left hand side of the above equation is in the range by construction. Hence, we find that
\par\
$f(z_j,z_n)=z_j^m \overline{z_j}^k \overline{z_n}^{l+1}= L_{jn}(\frac{1}{i(k+1)} z_j^m \overline{z_j}^{k+1} \overline{z_n}^{l+1})
 - L_{jn}(g) \in  L_{jn}(\mathcal{P})$, by the linearity of the operator.
\par\
Therefore, by the principle of induction and our preliminary arguments, every polynomial is in the range of $L_{jn}$, and as $j=1,...,n-1$ was taken arbitrarily, our desired result is proven.
\par\
$\textbf{\emph{QED}}$
\par\ \par\
From Rudin's book (Theorem 18.1.9 in [6]) we may also assert the following fact:
\par\ \par\
$\textbf{\emph{\underline{Remark:}}}$ If a (sufficiently) smooth function $u: \mathbb{H}^{2n-1} \rightarrow \mathbb{C}$ satisfies the Cauchy-Riemann equations on $\mathbb{H}^{2n-1}$, i.e. if $L_{jn}(u) = 0$ (for all $j=1,...,n-1$) on $\mathbb{H}^{2n-1}$, then $u$ admits a holomorphic extension $U: \mathcal{U} \rightarrow \mathbb{C}$, where $\mathcal{U} = \{\rho < 0\} = \{Im(z_n) \geq \Sigma_{j=1}^{n-1} \abs{z_j}^2 \} \subset \mathbb{C}^n$. More precisely, $U$ is holomorphic in $Int(\mathcal{U})$ and continuous to the boundary surface $\mathbb{H}^{2n-1}$, with $U|_{\mathbb{H}^{2n-1}} = u$ .
\par\ \par\
Let us now turn our attention to our problem of interest, namely the structure of the set of complex tangents to embeddings of $\mathbb{H}^{2n-1} \hookrightarrow \mathbb{C}^{2n-1}$. As we discussed above, the Heisenberg group $\mathbb{H}^{2n-1}$ arises naturally as a hypersurface of $\mathbb{C}^n$. In our paper [2], we considered embeddings of $\mathbb{H}^3 \hookrightarrow \mathbb{C}^3$ as graphs of maps $f: \mathbb{C}^2 \rightarrow \mathbb{C}$, i.e. $graph(f|_{\mathbb{H}^3}): \mathbb{H}^3 \hookrightarrow \mathbb{C}^3$. In the higher dimensional case, we need to take "several" graphs to achieve the desired embedding; in particular choose $n>2$ and let $f_1,...,f_{n-1}: \mathbb{C}^n \rightarrow \mathbb{C}$ be (sufficiently) smooth maps. Consider $F=(f_1, ... , f_{n-1}): \mathbb{C}^n \rightarrow \mathbb{C}^{n-1}$. Then $graph(F|_{\mathbb{H}^{2n-1}}): \mathbb{H}^{2n-1} \hookrightarrow \mathbb{C}^{2n-1}$ is an embedding of the desired type.
\par\ \par\
Now, we are interested in computing the complex tangents to such embeddings. We refer the reader to Webster in [7] for more on the following general formulation.
\par\
Let $\mathcal{M} \subset \mathbb{C}^n$ be a real $n$-manifold in $\mathbb{C}^n$, given by equations: $\mathcal{M}= \{r_1=0, ..., r_n=0\}$, where $r_i:\mathbb{C}^n \rightarrow \mathbb{R}$ are smooth and $dr_1 \wedge ... \wedge dr_n \neq 0$ (nowhere zero on $\mathcal{M}$). Let $\partial r_1 \wedge ... \wedge \partial r_n = B dz_1 \wedge ... \wedge dz_n$; hence $B:\mathcal{M} \rightarrow \mathbb{C}$ is a smooth map, in fact $B$ is the determinant of the $(nxn)$-matrix: $B(x) = det(\frac{\partial(r_1, ... , r_n)}{\partial(z_1, ... , z_{n})})$. Then a point $x \in \mathcal{M}$ is complex tangent if and only if $B(x)=0$.
\par\ \par\
Let's now apply this formulation for the situation at hand. In fact, our result for the 3-dimensional Heisenberg group in [2] extends in the analogous manner to all Heisenberg groups, as we state (and prove) in the following theorem:
\par\ 
\begin{theorem}: Let $M= \mathbb{H}^{2n-1}$ for $n \geq 2$ be any Heisenberg group. Then every (real) non-empty algebraic set that is the solution of one or two real polynomial equations in  $\mathbb{H}^{2n-1}$ may arise precisely as the set of complex tangents to some (algebraic) embedding  $\mathbb{H}^{2n-1} \hookrightarrow \mathbb{C}^{2n-1}$. In particular, every Heisenberg group admits a totally real embedding.
 \end{theorem}
\par\ \par\
$\textbf{\emph{\underline{Proof:}}}$ The situation for $n=2$ was proven in [2]. We will proceed to reduce our problem for $n>2$ to two variables and use the proof for the case $n=2$ to establish the result in all higher dimensions.
\par\ \par\
Let $M^{2n-1} = \{\rho=0\} \subset \mathbb{C}^n$ be a real hypersurface. Consider $\widetilde{M} = graph(F|_M) \subset \mathbb{C}^{2n-1}$, where $F=(f_1,...,f_{n-1}):\mathbb{C}^n \rightarrow \mathbb{C}^{n-1}$ is sufficiently smooth (at least continuously differentiable). Giving $\mathbb{C}^{2n-1}$ (holomorphic) coordinates $(z_1,...,z_n,\zeta_1, ..., \zeta_{n-1})$, we may write: $\widetilde{M} = \{R_1=0,..., R_{n-1} =0,R_n=0\}$, where: $R_n=\rho, R_1=\zeta_1-f_1, ..., R_{n-1}=\zeta_{n-1} - f_{n-1}$. Then we may describe $\widetilde{M} = \{Re(R_1)=0, Im(R_1) = 0, ..., Re(R_{n-1})=0, Im(R_{n-1}) = 0, R_{n-1}=0\}$ as the solution set of $n$ real equations, as prescribed in the formulation we put forth before. We are interested in the computation of the determinant: $B(x) = det(\frac{\partial(r_1, ... , r_{2n-1})}{\partial(z_1, ... , z_{n}, \zeta_1,..., \zeta_{n-1})})$, where $r_1= Re(R_1), r_2 = Im(R_1),...$ etc.  A point $x \in \mathcal{M}$ is complex tangent if and only if $B(x)=0$.
\par\ \par\
Let us choose: $f_2=\overline{z_1}, ...., f_{n-1} = \overline{z_{n-2}}$ and leave $f_1 = f$ arbitrary.
\par\
Then the $(2n-1) \times (2n-1)$ matrix simplifies and we proceed with our computation for $B$:
\par\ \par\
$B = det \left( {\begin{array}{cccccccccc}
   -f_{z_1} & -\overline{f}_{z_1} & 0 & -1 & 0 & 0 &... & 0 & 0 & 2\overline{z_1}\\
   -f_{z_2} & -\overline{f}_{z_2} & 0 & 0 & 0 & -1 & ... & 0 & 0 & 2 \overline{z_2}\\
   ... & ... & ... & ... & ... & ... & ... & ... & ... & ... \\
   -f_{z_{n}} & -\overline{f}_{z_{n}} & 0 & 0 & 0 & 0  &... & 0 & 0 & i \\
   1 & 0 & 0 & 0 & 0 & 0 & ... & 0 & 0 & 0 \\
  0 & 0 & 1 & 0 & 0 & 0 & ... & 0 & 0 & 0 \\
   ... & ... & ... & ... & ... & ... & ... & ... & ... & ... \\
  0 & 0 & 0 & 0 & 0 & 0 & ... & 1 & 0 & 0 \\
   \end{array} } \right) $
   \par\ \par\ \par\
Now, applying the determinant formula to the last $n-1$ rows, our computation simplifies considerably:
\par\ \par\
$=  det \left( {\begin{array}{ccccccc}
    -\overline{f}_{z_1} & -1 & 0 & ... &  0 & 0 & 2\overline{z_1}\\
    -\overline{f}_{z_2} & 0 & -1 & ... & 0 & 0 & 2 \overline{z_2}\\
    ... & ... & ... & ... & ... & ... & ... \\
    -\overline{f}_{z_{n-2}} & 0 & 0  &... & 0 & -1 & 2 \overline{z_{n-2}} \\
    -\overline{f}_{z_{n-1}} & 0 & 0  &... & 0 & 0 & 2 \overline{z_{n-1}} \\
    -\overline{f}_{z_{n}} & 0 & 0 &... & 0 & 0 & i \\
   \end{array} } \right) $.
\par\ \par\ \par\
Inserting the last column as the second column, we get
\par\ \par\
$= (-1)^n det \left( {\begin{array}{ccccccc}
    -\overline{f}_{z_1} & 2\overline{z_1} & -1 & 0 & ... &  0 & 0 \\
    -\overline{f}_{z_2} & 2 \overline{z_2} & 0 & -1 & ... & 0 & 0 \\
    ... & ... & ... & ... & ... & ... & ... \\
    -\overline{f}_{z_{n-2}}  & 2 \overline{z_{n-2}} & 0 & 0  &... & 0 & -1 \\
    -\overline{f}_{z_{n-1}} & 2 \overline{z_{n-1}}  & 0 & 0  &... & 0 & 0 \\
    -\overline{f}_{z_{n}} & i & 0 & 0 &... & 0 & 0  \\
   \end{array} } \right) $
\par\ \par\ \par\
$= \pm ( 2 \overline{f}_{z_{n}} \overline{z_{n-1}} - i \overline{f}_{z_{n-1}})$
\par\ \par\
Now, as we are only interested in the zeros of $B$, we may just as well consider $W = \pm \overline{B}$.
\par\
We now obtain: $W(f) = 2 z_{n-1} f_{\overline{z_n}} + i f_{\overline{z_{n-1}}}$.
\par\
But $W(f) = L_{n-1,n} (f)$, which is one of the "standard" tangential CR-operators to  $\mathbb{H}^{2n-1}$. We thus know $W$ is onto from Lemma 1, and for any complex polynomial $p(z_1, ..., z_n)$ there exists some polynomial f so that $W(f) = p$. Again, let $\mathcal{A} \subset \mathbb{H}^{2n-1}$ be a real algebraic set of codimension 1 or 2, that is: $\mathcal{A} = \{q_1 =0\} \bigcap \{q_2 = 0\}  \bigcap \mathbb{H}^{2n-1}$ for some polynomials $q_1, q_2: \mathbb{R}^{2n} \rightarrow \mathbb{R}$ (potentially equal). Then, taking $p=q_1 + i q_2$ and changing variables accordingly, $p(z_1, ..., z_n)$ is a complex polynomial. Hence, there exists a polynomial $f$ so that $W(f) = p$. Taking the embedding: $E = graph(F|_{\mathbb{H}^{2n-1}}): \mathbb{H}^{2n-1} \hookrightarrow \mathbb{C}^{2n-1}$, $F= (f, \overline{z_1}, ..., \overline{z_{n-2}})$ we see from the above that the complex tangents to the the embedding E will be precisely the zeros of $W(f) = p$ on $\mathbb{H}^{2n-1}$. But $\{p=0\} \bigcap \mathbb{H}^{2n-1} = \{q_1 =0\} \bigcap \{q_2 = 0\} \bigcap \mathbb{H}^{2n-1} = \mathcal{A}$! Hence, every algebraic subset given by one or two real polynomial equations in $\mathbb{H}^{2n-1}$ may arise precisely as the set of complex tangents to some (polynomial) embedding.
\par\ \par\
Furthermore, taking $f=\overline{z_{n-1}}$, the embedding given by the graph of:  $F= (\overline{z_{n-1}}, \overline{z_1}, ..., \overline{z_{n-2}})$ must necessarily be totally real, as $W(\overline{z_{n-1}}) = i \neq 0$.
\par\
$\textbf{\emph{QED}}$
\par\ \par\
\section*{2. Results for Odd-Dimensional Spheres}

We may further extend our results to odd dimensional spheres in an analogous manner to what we did for $S^3$ (see [2]). For any fixed $n \geq 2$, there exists a natural biholomorphism between the sphere $S^{2n-1}$ minus its "north pole" and the Heisenberg group $\mathbb{H}^{2n-1}$, namely let $\varphi: \mathbb{H}^{2n-1} \rightarrow S^{2n-1} \setminus \{(0,..., 0, 1)\}$ be the map:
\par\
$\varphi(z_1,...,z_n) = (\frac{2z_1}{z_n+i}, ..., \frac{2z_{n-1}}{z_n+i}, \frac{z_n-i}{z_n+i})$.
\par\
And let $\psi: S^{2n-1} \setminus \{(0, ..., 0, 1)\} \rightarrow \mathbb{H}^{2n-1}$ be its inverse. So:
\par\
$ \psi(z_1, ..., z_n) = (\frac{iz_1}{1-z_n}, ..., \frac{iz_{n-1}}{1-z_n}, i\frac{1+z_n}{1-z_n})$.
\par\ \par\
We may then use these biholomorphisms to extend our results on the higher dimensional Heisenberg groups to all odd dimensional spheres, analogous to our constructs in dimension 3 ($n=2$). In particular we show that every odd dimensional sphere admits an embedding $S^{2n-1} \hookrightarrow \mathbb{C}^{2n-1}$ which assumes its complex tangents along any prescribed algebraic set. More precisely:
\par\ 
\begin{theorem} Let $n \geq 2$ and consider the sphere $S^{2n-1} \subset \mathbb{C}^n$. Then every (real) non-empty algebraic set given by one or two polynomial equations in  $S^{2n-1}$ may arise precisely as the set of complex tangents to some $\mathcal{C}^k$-embedding  $S^{2n-1} \hookrightarrow \mathbb{C}^{2n-1}$, for any prescribed $k \in \mathbb{N}$.
 \end{theorem}
\par\ \par\
$\textbf{\emph{\underline{Proof:}}}$ The case for the 3-sphere (n=2) was already proved in our paper [2], so assume $n > 2$.
\par\
Let $\mathcal{A} \subset S^{2n-1}$ be a (non-empty) real algebraic set given by one or two polynomial equations, that is  $\mathcal{A} = \{q_1 = 0\} \bigcap \{q_2 = 0\} \bigcap S^{2n-1}$, where $q_1, q_2 : \mathbb{C}^3 \rightarrow \mathbb{R}$ are polynomials (possibly equal). Then again we may write $p = q_1 + i q_2$ and change coordinates to get a single complex polynomial in $(z_1, z_2, z_3)$ variables (and their anti-holomorphic counterparts). So $\mathcal{A} = \{p = 0\} \bigcap S^{2n-1}$. Furthermore, via a (polynomial) rotation, we may insure that $\mathcal{A}$ contains the north pole $(0, ..., 0, 1)$.
\par\ \par\
Now, consider the set $\psi (\mathcal{A} \setminus \{(0, ..., 0, 1)\}) \subset \mathbb{H}^{2n-1}$. This is an algebraic set; in fact we demonstrated in section 3 of our paper [1] that birational maps preserve algebraic sets, and is direct that it will be given by the same number of real polynomial equations (one or two) as $\mathcal{A}$ is defined a subset of $S^{2n-1}$.
\par\
By the Theorem 2 above, there exists some polynomial $g: \mathbb{C}^n \rightarrow \mathbb{C}$ so that the complex tangents to embedding given by the graph of the function: $G = (g, \overline{z_1}, ..., \overline{z_{n-2}})$ restricted to $\mathbb{H}^{2n-1}$ form precisely the given algebraic set $\psi (\mathcal{A} \setminus \{(0, ..., 0, 1)\})$. We can explicitly write this embedding. Let
\par\ \par\
$\mathcal{M}_g = \{(z_1, ..., z_n, g(z_1, ..., z_n), \overline{z_1}, ..., \overline{z_{n-2}}) | (z_1, ..., z_n) \in \mathbb{H}^{2n-1}\} \subset \mathbb{C}^{2n-1}$
 \par\
 be the graph indicated above. Then the set of complex tangents to $\mathcal{M}_g$ is $\psi (\mathcal{A} \setminus \{(0, ..., 0, 1)\})$.
\par\ \par\
Consider now the mapping: $g \circ \psi: S^{2n-1} \setminus \{(0,...,0,1)\}\rightarrow \mathbb{H}^{2n-1} \rightarrow \mathbb{C}$. This mapping is smooth (on $S^{2n-1}$) and in fact takes the form of the rational function:
\par\
$(g \circ \psi) (z_1, ..., z_n) = \frac{s(z_1,...,z_n)}{(1-z_n)^l (1-\overline{z_n})^l}$
\par\
where $l \in \mathbb{N}$ is the degree of the polynomial $g$.
\par\
Now, consider the $graph((G \circ \psi)|_{S^{2n-1} \setminus \{(0, ..., 0, 1)\}}) \subset \mathbb{C}^{2n-1}$, which we write:
\par\ \par\
$\mathcal{S}_{g \circ \psi} = \{(w_1,...,w_n, g(\psi(w_1,...,w_n)), \overline{\psi_1 (w_1,...,w_n)},..., \overline{\psi_{n-2} (w_1, ..., w_n)}) |$
\par\
$(w_1,...., w_n) \in S^{2n-1} \setminus \{0,...,0,1\}\} \subset \mathbb{C}^{2n-1}$.
\par\
 Note that $\psi_j (w_1,...,w_n) = i\frac{w_j}{1-w_n}$ is the $j^{th}$-coordinate of the biholomorphism $\psi$.
\par\ \par\
The spaces $\mathcal{M}_g, \mathcal{S}_{g \circ \psi} \subset \mathbb{C}^{2n-1}$ are then biholomorphic, via $\Lambda: \mathcal{M}_g \rightarrow \mathcal{S}_{g \circ \psi}$, acting by the biholomorphism $\varphi:\mathbb{H}^{2n-1} \rightarrow S^{2n-1}$ on the first $n$ coordinates, and by identity on the remaining coordinates.
\par\
Hence, $\Lambda$ preserves the complex tangents of the respective spaces, and as the complex tangents of the (graphical) embedding: $\mathbb{H}^{2n-1} \cong \mathcal{M}_g \subset  \mathbb{C}^{2n-1}$ form the set: $\psi (\mathcal{A} \setminus \{(0, ..., 0, 1)\})$, the complex tangents of the embedding:
 \par\
 $S^{2n-1} \{(0,...,0,1)\} \cong \mathcal{S}_{g \circ \psi} \subset  \mathbb{C}^{2n-1}$
 \par\
 will form the set: $\varphi(\psi (\mathcal{A} \setminus \{(0, ..., 0, 1)\})) = \mathcal{A} \setminus \{(0, ..., 0, 1)\} \subset S^{2n-1}$.
\par\ \par\
Hence, $\mathcal{S}_{g \circ \psi}$ gives us an embedding of the sphere $S^{2n-1}$ minus its north pole whose complex tangents are precisely our ordained algebraic set $\mathcal{A} \setminus \{(0, ..., 0, 1)\}$ (minus the north pole). Unfortunately, the manifold:
 \par\
 $\mathcal{S}_{g \circ \psi} = graph((G \circ \psi)|_{S^{2n-1} \setminus \{(0, ..., 0, 1)\}}) \subset \mathbb{C}^{2n-1}$ assumes a non-removable singularity at the north pole; in fact the mapping diverges there. Consider, however the function: (for any $r \in \mathbb{N}$)
\par\ \par\
$ \mathcal{G}_r (z_1, ..., z_n)= (1-z_n)^{2l+r} \cdot (G \circ \psi)(z_1, ...., z_n)$
\par\
$= (\frac{(1-z_n)^{l+r}}{(1-\overline{z_n})^l} s(z_1, ..., z_n), \frac{-i\overline{z_1}(1-z_n)^{2l+r}}{1-\overline{z_n}}, ..., \frac{-i\overline{z_{n-2}}(1-z_n)^{2l+r}}{1-\overline{z_n}})$,
\par\
in which we multiply in every coordinate by the holomorphic factor $(1-z_n)^{2l+r}$.
\par\ \par\
Then $\Upsilon: graph(\mathcal{G}_r |_{S^{2n-1} \setminus \{(0, ..., 0, 1)\}}) \cong graph((G \circ \psi)|_{S^{2n-1} \setminus \{(0, ..., 0, 1)\}}) \subset \mathbb{C}^{2n-1}$ and hence the spaces are biholomorphic. Note that this mapping is given as the identity in first $n$ coordinates and multiplication by the holomorphic factor  $(1-z_n)^{2l+r}$ in each of the other coordinates. Clearly the map is one-to-one and onto by construction, and holomorphic with holomorphic inverse.
\par\ \par\
Hence, the set: $ graph(\mathcal{G}_r |_{S^{2n-1} \setminus \{(0, ..., 0, 1)\}}) \subset \mathbb{C}^{2n-1}$ has its complex tangents along $\mathcal{A} \setminus \{(0, ..., 0, 1)\}$ as well. Furthermore, the map $\mathcal{G}_r$ can be defined at the north pole $(0,...,0,1)$ and takes the value zero, with $\mathcal{G}_r$ of class $\mathcal{C}^{r-1}$ as a map on the sphere.
\par\ \par\
Therefore, the embedding: $S^{2n-1} \cong \widetilde{\mathcal{S}}_r = graph(\mathcal{G}_r |_{S^{2n-1}}) \subset \mathbb{C}^{2n-1}$ is an embedding of the sphere continuously differentiable of class $\mathcal{C}^{r-1}$ whose complex tangents away from the north pole assume the algebraic set $\mathcal{A} \setminus \{(0, ..., 0, 1)\}$. Now, it remains to resolve the issue at the north pole.
\par\ \par\
Recall from our previous constructions that the complex tangents to such a graphical embedding will be given by the zeros of the determinant of the $(n \times n)$-matrix with the first $n-1$ columns being the (holomorphic) gradients of the conjugates of the component functions and the last column being the vector: $(z_1, ..., z_n)$. Now, applying the gradient formula to the last column, we get a sum of determinants of $(n-1)\times(n-1)$-matrices whose entries are all (single) holomorphic derivatives of the component functions. In our case, by the construction of the component functions for $\mathcal{G}_r$, out of each of the entries we can factor out a term: $\frac{(1-\overline{z_n})^{l+r-1}}{(1-z_n)^{l+1}}$ (after applying the conjugation). Assuming without loss of generality that $r>2$, we find that each term must evaluate to be zero at the north pole $(0,...,0,1)$. As such, the determinant of the relevant matrix must be zero at the north pole.
\par\ \par\
Hence, by construction, we find that the point $(0,...,0,1)$ is also complex tangent. By our results above, we conclude that the (graphical) embedding:
\par\
$S^{2n-1} \cong \widetilde{\mathcal{S}}_r \subset \mathbb{C}^{2n-1}$ will be an embedding of class $\mathcal{C}^{r-1}$ assuming its complex tangents exactly along the algebraic set $\mathcal{A}$. As $r>2$ and $\mathcal{A} \subset S^{2n-1}$ were chosen arbitrarily, our claim is proven.
\par\
$\textbf{\emph{QED}}$
\par\ \par\
We recall Gromov's result that the only spheres $S^k$ that admit totally real embeddings into $\mathbb{C}^k$ are the spheres $S^1, S^3$ (the case $S^1$ being trivial). Our following corollary show that the "minimal configuration" for complex tangents is possible for all odd-dimensional spheres; more precisely:
\par\
\begin{corollary}: Let $n \geq 2$ and consider the sphere $S^{2n-1} \subset \mathbb{C}^n$. Then $S^{2n-1}$ admits an embedding $S^{2n-1} \hookrightarrow \mathbb{C}^{2n-1}$ who assumes exactly one complex tangent point (at the north pole). Further, this embedding may be assumed to be of class $\mathcal{C}^r$ for any $r \in \mathbb{N}$.
 \end{corollary}
\par\
$\textbf{\emph{\underline{Proof:}}}$ Take $g=\overline{z_{n-1}}$ in the above proof. Then the embedding given by the graph of:
\par\
$ \mathcal{G}_r (z_1, ..., z_n) = (\frac{-i\overline{z_{n-1}}(1-z_n)^{2+r}}{1-\overline{z_n}}, \frac{-i\overline{z_1}(1-z_n)^{2+r}}{1-\overline{z_n}}, ..., \frac{-i\overline{z_{n-2}}(1-z_n)^{2+r}}{1-\overline{z_n}})$.
\par\
will give an embedding of the desired type, i.e. complex tangent only at the north pole $(0, ..., 0, 1)$ and of class $\mathcal{C}^r$ at the north pole, and smooth otherwise.
\par\
$\textbf{\emph{QED}}$
\par\ \par\
$\textbf{\emph{\underline{Remark:}}}$ We may further extend our results to some "Heisenberg-like" spaces. In particular, consider first in $\mathbb{C}^2$ a hypersurface $M=\{\rho =0\}$, where:
\par\
$\rho (z,w) = \alpha Re(z) + \beta Im(z) + h(w)$, where $\alpha, \beta \in \mathbb{R}$ and $h$ is a complex polynomial (in $w, \overline{w}$) with: $\overline{h} = h$.
\par\
Note such hypersurfaces in a sense take the form of a cylinder as hypersurfaces in $\mathbb{R}^4$ ($=\mathbb{C}^2$).
\par\ \par\
It is direct to see that our arguments given for the Heisenberg group $\mathbb{H}^3$ can be readily extended to the $M$, and we obtain that every complex polynomial will be in the range of the (single) tangential CR-operator to $M$, which we denoted by $\mathcal{L}_M$. As such, every real algebraic set given by one or two real polynomial equations (defined in $\mathbb{R}^4$) intersected with $M$ may arise as the set of complex tangents to some embedding into $\mathbb{C}^3$. In particular, all such "trough-like" surfaces admit totally real embeddings, which can be taken to be polynomial embeddings if either $\alpha$ or $\beta$ is not equal zero.
\par\ \par\
Using our constructions given previously in this section, we may further extend our results for higher dimensional analogues. In particular consider any "trough-like" hypersurface of the form: $M=\{\rho=0\} \subset \mathbb{C}^n$, given by $\rho(z_1,...,z_n)= \alpha Re(z_1) + \beta Im(z_1) + h(z_2,...,z_n)$, where $\alpha, \beta \in \mathbb{R}$ and $h$ is a complex polynomial (in $z_j, \overline{z_j}$, $2 \leq j \leq n$) with: $\overline{h} = h$. Then any algebraic set given by one or two real polynomial equations in $M$ may arise precisely as the set of complex tangents to some embedding: $M \hookrightarrow \mathbb{C}^{2n-1}$. Again, it follows that all such hypersurfaces will admit totally real embeddings.

\par\ \par\

\end{document}